\documentclass[12pt,leqno]{amsart}
\usepackage{amsmath}
\usepackage{amssymb}

\def\underset#1#2{{\mathrel{\mathop {{}_{} {#2}}\limits_{{#1}_{}}}}}
\def\upplim_#1{\underset{#1}{\overline\lim}\;}
\def\lowlim_#1{\underset{#1}{\underline\lim}\;}

\def\le{\leqslant}
\newtheorem{corollary}[equation]{Corollary}

\newtheorem{claim}[equation]{Claim}

\newtheorem{theorem}[equation]{Theorem}

\newcommand{\C}{{\mathbf{C}}}

\renewcommand{\P}{{\mathbf{P}}}

\newcommand{\supp}{\mathrm{Supp}\,}

\numberwithin{equation}{section}

\begin{document}
\noindent

\title[Algebraic dependences of meromorphic mappings]{Algebraic dependences and uniqueness problem of meromorphic mappings sharing moving hyperplanes without counting multiplicities }

\author{Le Ngoc Quynh}
\address{Faculty of Education, An Giang University, 18 Ung Van Khiem, Dong Xuyen, Long Xuyen, An Giang, Vietnam}
\email{nquynh1511@gmail.com}
\maketitle
\begin{abstract}
This article deals with the multiple values and algebraic dependences problem of meromorphic mappings sharing moving hyperplanes in projective space. We give some algebraic dependences theorems for meromorphic mappings sharing moving hyperplanes without counting multiplicity, where all zeros with multiplicities more than a certain number are omitted. Basing on these results, some unicity theorems regardless of multiplicity for meromorphic mappings in several complex variables are given. These results are extensions and strong improvements of some recent results.
\end{abstract} 

\def\thefootnote{\empty}

\section{Introduction}

The theory on algebraic dependences of meromorphic mappings in several complex variables into the complex projective spaces for fixed targets was first studied by W. Stoll \cite{St1}. Later, M. Ru \cite{R} generalized W. Stoll's result to the case of holomorphic curves into the complex projective spaces sharing moving hyperplanes. Recently, by using the new second main theorem given by Thai-Quang \cite{TQ}, P. D. Thoan, P. V. Duc and S. D. Quang \cite{PP,PPS,Q2} gave some  improvements of the results of W. Stoll and M. Ru. In order to state some of  their result, we first recall the following.

We call a meromorphic mapping of $\C^n$ into $\P^N(\C)^*$ a moving hyperplane in $\P^N(\C)$.
 Let $a_1,\dots,a_q$ $(q \geq N+1)$ be $q$ moving hyperplanes with reduced representations $a_j = (a_{j0}: \dots : a_{jN})\ (1\le j \le q).$ We say that $a_1,\dots,a_q$ are in general position if $\det (a_{j_kl}) \not \equiv 0$ for any $1\le j_0<j_1<...<j_N\le q.$

Let $f_i: \C^n\rightarrow \P^N(\C) \ (1\leqslant i \leqslant \lambda)$ be meromorphic mappings with reduced representations
 $f_i:=(f_{i0}:\cdots :f_{iN}).$ Let $g_j: \C^n\rightarrow \P^N(\C)^* \ (0\leqslant j\leqslant q-1)$ be meromorphic mappings in
general position with reduced representations $g_j:=(g_{j0}:\cdots :g_{jN}).$ Put $(f_i,g_j):= \sum_{s=0}^Nf_{is}g_{js}\ne 0$ for each  
$1\le i\le\lambda ,\ 0\le j\le q-1$ and assume that $\min\{1,\nu^0_{(f_1,g_j),\le k_j}\}=\cdots =\min\{1,\nu^0_{(f_\lambda,g_j),\le k_j}\}$. Put $A_j=\supp (\nu^0_{(f_1,g_j),\le k_j}).$ Assume that each $A_j$ has an irreducible decomposition as follows $A_j=\bigcup_{s=1}^{t_j}A_{js}$.
Set $A=\bigcup_{A_{js}\not\equiv A_{j's'}}\{A_{js} \cap A_{j's'}\}$ with $1\le s\le t_j,1\le s'\le t_{j'}, 0\le j, j' \le q-1$.

Denote by $T[N+1,q]$ the set of all injective maps from $\{1,\cdots,N+1\}$ to $\{0,\cdots ,q-1\}.$ 
For  each $z\in \C^n\setminus\{\bigcup_{\beta\in T[N+1,q]}\{z|g_{\beta (1)}(z)\wedge\cdots\wedge g_{\beta (N+1)}(z)=0\} \cup A\cup\bigcup_{i=1}^{\lambda}I(f_i)\},$ we define $\rho (z)=\sharp\{j|z\in A_j\}$. Then $\rho (z)\le N.$ For any positive number $r>0,$ define $\rho (r)=\sup\{\rho (z)| |z|\le r\},$ where the supremum is taken over all
$z\in \C^n\setminus\{\bigcup_{\beta\in T[N+1,q]}\{z|g_{\beta (1)}(z)\wedge\cdots\wedge g_{\beta (N+1)}(z)=0\} \cup A\cup\bigcup_{i=1}^{\lambda}I(f_i)\}.$ 
Then $\rho (r)$ is a decreasing function. Let $$ d:=\lim_{r\rightarrow +\infty}\rho (r).$$
Then $d\le N.$  If for each $i \ne j, \ \dim\{A_i\cap A_j\}\le n-2,$ then $d=1.$

When all $k_j=+\infty$, in \cite{PP}, P. V. Duc and P. D. Thoan proved the following.
\vskip0.2cm
\noindent
{\bf Theorem A} (see \cite[Theorem 1]{PP}).\ {\it Let $f_1,\cdots ,f_{\lambda}: \C^n\rightarrow \P^N(\C)$ be non-constant meromorphic mappings. Let $g_i:\C^n \rightarrow \P^N(\C)\ (0\le i\le q-1)$
be moving targets located in general position and $T(r,g_i)=o(\max_{1\le j\le\lambda}T(r,f_j))\ (0\le i\le q-1).$ Assume that $(f_i,g_j)\not\equiv 0$
for $1\le i\le\lambda ,\ 0\le j\le q-1$ and $A_j:=(f_1,g_j)^{-1}\{0\}=\cdots =(f_{\lambda},g_j)^{-1}\{0\}$ for each $0\le j\le q-1.$ Denote $\mathcal A=\cup_{j=0}^{q-1}A_j$. Let $l, \ 2\le l\le \lambda,$ be an integer such that for any increasing sequence  $1\le j_1<\cdots <j_l\le\lambda,\ f_{j_1}(z)\wedge\cdots\wedge f_{j_l}(z)=0$ for every point $z\in\mathcal A.$ If $q>\dfrac{dN(2N+1)\lambda}{\lambda -l+1},$ then $f_1,\cdots ,f_{\lambda}$ are algebraically over $\C,$ i.e.  $f_1\wedge\cdots\wedge f_{\lambda}\equiv 0$.}

\vskip0.2cm
Furthermore, in the case of $d=1$, P. D. Thoan, P. V. Duc and S. D. Quang \cite{PPS} proved an better algebraic dependences theorem as follows.

\vskip0.2cm
\noindent
{\bf Theorem B}\ (see \cite[Theorem 1]{PPS}). \ {\it  Let $f_1,\cdots ,f_{\lambda}: \C^n\rightarrow \P^N(\C)$ be non-constant meromorphic mappings. 
Let $g_i:\C^n \rightarrow \P^N(\C)^* \ (0\le i\le q-1)$ be moving hyperplanes in general position such that $T(r,g_i)=o(\max_{1\le j\le\lambda}T(r,f_j))\ (0\le i\le q-1)$ and $(f_i,g_j)\not\equiv 0$ for $1\le i\le\lambda ,\ 0\le j\le q-1.$ Assume that the following conditions are satisfied.
\begin{itemize}
\item[(a)] $\min\{1,\nu_{(f_1,g_j)}\}=\cdots =\min\{1,\nu_{(f_{\lambda},g_j)}\}$ for each $0\le j\le q-1,$ 
\item[(b)] $\dim\{z|(f_1,g_i)(z)=(f_1,g_j)(z)=0\}\le n-2$ for each $0\le i< j\le q-1,$
\item[(c)] there exists an integer number $l, \ 2\le l\le \lambda,$ such that for any increasing sequence  $1\le j_1<\cdots <j_l\le\lambda,\ f_{j_1}(z)\wedge\cdots\wedge f_{j_l}(z)=0$ for every point $z\in \bigcup_{i=0}^{q-1}(f_1,g_i)^{-1}\{0\}.$
\end{itemize}
If $q>\frac{N(2N+1)\lambda -(N-1)(\lambda -1)}{\lambda -l+1},$ then $f_1\wedge\cdots\wedge f_{\lambda}\equiv 0$.}

\vskip0.2cm 
The above results are the best results on the algebraic dependences of meromorphic mappings sharing moving hyperplanes available at the present. However, in our opinion, they are still weak. Also in the above results, all intersecting points of the mappings and the moving hyperplanes are considered. Actually, there are many authors consider the multiple values for meromorphic mappings sharing hyperplanes, i.e., consider only the intersecting points of the mappings $f_i$ and the hyperplanes $g_j$ with the multiplicity not exceed a certain number $k_j<+\infty$. For example, in 2010, T. B. Cao and H. X. Yi gave some uniqueness theorems for meromorphic mappings sharing fixed hyperplanes where all intersecting points more than a certain number are omitted (see \cite[Theorems 1.4 and 1.5]{CY}). Recently, H. H. Giang consider the multiple values and uniqueness problems for the mappings sharing moving hyperplanes (see \cite[Theorems 1.1 and 1.3]{H}).

Our first purpose in this paper is to generalize and improve Theorems A and B by considering the multiple values problem and reducing the number of hyperplanes. Namely, we will prove the following.

\begin{theorem}\label{1.1}
Let $f_1,\cdots ,f_{\lambda}: \C^n\rightarrow \P^N(\C)$ be non-constant meromorphic mappings. Let $\{g_j\}_{j=0}^{q-1}$ be moving hyperplanes of $\P^N(\C)$ in general position satisfying $T(r,g_j)=o(\max_{1\le i\le\lambda}T(r,f_i))\ (0\le j\le q-1).$ Let $k_j \ (0\le j\le q-1)$ be positive integers or $+\infty.$ Assume that $(f_i,g_j)\not\equiv 0$ for $1\le i\le\lambda ,\ 0\le j\le q-1$ and $A_j:=\supp \nu^0_{(f_1,g_j),\le k_j} = \cdots = \supp \nu^0_{(f_\lambda,g_j),\le k_j}$ for each $0\le j\le q-1.$ Denote $\mathcal A=\bigcup_{j=0}^{q-1}A_j$. Let $l, \ 2\le l\le \lambda,$ be an integer such that for any increasing sequence  $1\le i_1<\cdots <i_l\le\lambda,\ f_{i_1}(z)\wedge\cdots\wedge f_{i_l}(z)=0$ for every point $z\in\mathcal A.$ 
If 
$$\sum\limits_{j=0}\limits^{q-1}\dfrac{1}{k_j} < \dfrac{q}{N(N+2)} - \dfrac{d\lambda}{\lambda-l+1},$$
then $f_1\wedge\cdots\wedge f_{\lambda}\equiv 0$.
\end{theorem}

In the above result, letting $k_j=+\infty \ (0\le j\le q-1)$, we get the conclusion of Theorem A with $q>\frac{d\lambda N(N+2)}{(\lambda - l+1)}.$ Also, letting $d=1$, we obtain the conclusion of Theorem B with $q>\frac{\lambda N(N+2)}{\lambda -l+1}$. Hence our result is improvement of Theorems A and B.

Let $\lambda =l=2$ and $k_j=+\infty \ (0\le j\le q-1)$, the above theorem implies the following unicity theorem.
\begin{corollary}\label{1.2}
Let $f_1, f_2:\C^n\rightarrow \P^N(\C)$ be two non-constant meromorphic mappings. Let $\{g_j\}_{j=0}^{q-1}$ be moving hyperplanes of $\P^N(\C)$ in general position satisfying $T(r,g_j)=o(\max_{1\le i\le 2}T(r,f_i))\ (0\le j\le q-1).$ Assume that the following conditions are satisfied.
\begin{itemize}
\item[(a)] $\min\{1,\nu_{(f_1,g_j)}^{0}\}=\min\{1,\nu_{(f_2,g_j)}^{0}\}$ for each $0\le j\le q-1,$
\item[(b)] $f_1(z)=f_2(z)$ for each $z\in\bigcup_{j=0}^{q-1}\supp \nu^0_{(f_1,g_j)}$,
\item[(c)] $q>2dN(N+2).$
\end{itemize}
Then $f_1\equiv f_2$.
\end{corollary}

If $d=1$, then from the above corollary, we get a uniqueness theorem for meromorphic mappings, which are not assumed to be nondegenerate, sharing $q>2N^2+4N$ moving hyperplanes in general position. 

\vskip0.2cm
Now, if we assume further that the linearly closures of the images of the mappings $f_i$ in Theorem \ref{1.1} have the same dimension then we will get a better result as follows.

\begin{theorem}\label{1.3}
 Let $f_1,\cdots ,f_{\lambda}: \C^n\rightarrow \P^N(\C)$ be non-constant meromorphic mappings. Let $\{g_j\}_{j=0}^{q-1}$ be $q$ meromorphic mappings $g_j:\C^n \rightarrow \P^N(\C)^* \ (0\le j\le q-1)$ in general position satisfying $T(r,g_j)=o(\max_{1\le i\le\lambda}T(r,f_i))\ (0\le j\le q-1).$ Let $k_j \ (0\le j\le q-1)$ be positive integers or $+\infty.$ Assume that $(f_i,g_j)\not\equiv 0$ for $1\le i\le\lambda ,\ 0\le j\le q-1,$ and the following conditions are satisfied.
\begin{itemize}
\item[(a)] $\min\{1,\nu^0_{(f_1,g_j)\le k_j}\}=\cdots =\min\{1,\nu^0_{(f_{\lambda},g_j)\le k_j}\}$ for each $0\le j\le q-1,$ 
\item[(b)] $\dim  {f_1}^{-1} (g_{j_1})\cap {f_1}^{-1}(g_{j_2}) \le n-2$ for each $0\le j_1< j_2\le q-1,$ 
\item[(c)] there exists an integer number $l, \ 2\le l\le \lambda,$ such that for any increasing sequence  $1\le i_1<\cdots <i_l\le\lambda,\ f_{i_1}(z)\wedge\cdots\wedge f_{i_l}(z)=0$ for every point $z\in \bigcup\limits_{j=0}^{q-1} {f_1}^{-1}(g_j).$ 
\end{itemize}
We assume further that $\mathrm{rank}_{\mathcal R\{g_j\}}f_1=\cdots =\mathrm{rank}_{\mathcal R\{g_j\}}f_\lambda =m+1$, where $m$ is a positive integer. If 
\begin{align}\tag{*}
\sum_{j=0}^{q-1}\dfrac{1}{k_j+1-m}<\dfrac{q}{m(2N-m+2)}-\dfrac{\lambda q}{q(\lambda -l+1)+\lambda (m-1)}
\end{align}
then $f_1\wedge\cdots\wedge f_{\lambda}\equiv 0$.
\end{theorem}

\noindent
\textbf{Remark 1:} i) When $k_0=\cdots =k_{q-1}=+\infty$, the condition (*) becomes 
$$q>\frac{\lambda (2Nm-m^2+m+1)}{(\lambda -l+1)}.$$
We see that this inequality is satisfied with $q>\frac{\lambda (N^2+N+1)}{(\lambda -l+1)}$, since its right hand side attains maximum at $m=N$. Hence in this case ($d=1$ and $k_0=\cdots =k_{q-1}=+\infty$), the conclusion of Theorem \ref{1.3} is better than that of Theorem \ref{1.1}. 

ii) Let $\lambda =l=2$ and $k_0=\cdots =k_{q-1}=+\infty$. We may show that if $f_1$ and $f_2$ satisfy the conditions (i)-(iii) of Theorem \ref{1.3} and $q>N(N+2)$ then $\mathrm{rank}_{\mathcal R\{g_j\}}f_1=\mathrm{rank}_{\mathcal R\{g_j\}}f_2$.

 Indeed, suppose that there are $a_0,...,a_N\in\mathcal R\{g_j\}$, not all zeros, satisfying $\sum_{0\le i\le N}a_if_{1i}\equiv 0$. Set $P=\sum_{0\le i\le N}a_if_{2i}$. Then $P=0$ on $\bigcup_{1\le j\le q-1}(f_1,g_j)^{-1}(0)$. If  $P\not\equiv 0$, then by using Remark 2 below, we have
\begin{align*}
T_{f_2}(r)&\ge N(r,\nu^0_P)+o(T_{f_2}(r))\ge\sum_{j=0}^{q-1}N^{[1]}(r,\nu^0_{(f_1,g_j)})+o(T_{f_2}(r))\\ 
&= \sum_{j=0}^{q-1}N^{[1]}(r,\nu^0_{(f_2,g_j)})+o(T_{f_2}(r))\ge \dfrac{q}{N(N+2)}T_{f_2}(r)+o(T_{f_2}(r)).
\end{align*}
Letting $r\longrightarrow +\infty$, we get $q\le N(N+2)$. This is a contradiction. Then we must have $P\equiv 0$. This implies that $\mathrm{rank}_{\mathcal R\{g_j\}}f_1\ge \mathrm{rank}_{\mathcal R\{g_j\}}f_2$.

Similarly, we have $\mathrm{rank}_{\mathcal R\{g_j\}}f_1\le \mathrm{rank}_{\mathcal R\{g_j\}}f_2$. Hence $\mathrm{rank}_{\mathcal R\{g_j\}}f_1= \mathrm{rank}_{\mathcal R\{g_j\}}f_2$

\vskip0.2cm
Then from Theorem \ref{1.3} and the above remarks, we get a uniqueness theorem as follows.
\begin{corollary}\label{1.4}
Let $f_1, f_2:\C^n\rightarrow \P^N(\C)$ be non - constant meromorphic mappings. Let $\{g_j\}_{j=0}^{q-1}$ be moving hyperplanes of $\P^N(\C)$ in general position satisfying $T(r,g_j)=o(\max_{1\le i\le 2}T(r,f_i))\ (0\le j\le q-1).$ Assume that $(f_i,g_j)\not\equiv 0$ for $1\le i\le 2,\ 0\le j\le q-1,$ and the following conditions are satisfied.
\begin{itemize}
\item[(a)] $\min\{1,\nu^0_{(f_1,g_j)}\}=\min\{1,\nu^0_{(f_2,g_j)}\}$ for each $0\le j\le q-1,$ 
\item[(b)] $\dim  {f_1}^{-1} (g_{j_1})\cap {f_1}^{-1}(g_{j_2}) \le n-2$ for each $0\le j_1< j_2\le q-1,$ 
\item[(c)] $f_1=f_2$ on $\bigcup\limits_{j=0}^{q-1}\supp \nu^0_{(f_1,g_j)}.$ 
\end{itemize}
If $q>2N^2+2N+2,$ then $f_1\equiv f_2.$
\end{corollary}

We would like to note that there are several results on the uniqueness problem of meromorphic mappings sharing moving hyperplanes regardless of multiplicity. For example, in 2007, with the same assumption of Corollary \ref{1.4} (and plus $N\ge 2$), Z. Chen, Y. Li and Q. Yan \cite{Chen} get the uniqueness theorem with $q\ge 4N^2+2N$. In 2013, P. D. Thoan, P. V. Duc and S. D. Quang improved the result of these authors to the case of $q\ge 4N^2+2$ (for any $N$). Therefore, our above uniqueness theorem is much stronger improvements of many previous results.

\noindent\textbf{Acknowledgements}: This paper was supported in part by a NAFOSTED grant of Vietnam.

\section{Basic notions and auxiliary results from Nevanlinna theory}

\noindent
{\bf 2.1.}\ We set $||z|| = \big(|z_1|^2 + \dots + |z_n|^2\big)^{1/2}$ for
$z = (z_1,\dots,z_n) \in \C^n$ and define
\begin{align*}
B(r) := \{ z \in \C^n : ||z|| < r\},\quad
S(r) := \{ z \in \C^n : ||z|| = r\}\ (0<r<\infty).
\end{align*}

Define 
$$v_{n-1}(z) := \big(dd^c ||z||^2\big)^{n-1}\quad \quad \text{and}$$
$$\sigma_n(z):= d^c \text{log}||z||^2 \land \big(dd^c \text{log}||z||^2\big)^{n-1}
\text{on} \quad \C^n \setminus \{0\}.$$

\noindent
{\bf 2.2.}\ For a divisor $\nu$ on $\C^n$, we denote by $N(r,\nu)$ the counting function of the divisor $\nu$ as usual in Nevanlinna theory (see \cite{NO}).

For a positive integer $M$ or $M= \infty$, we define the truncated divisors of $\nu$ by
$$\nu^{[M]}(z)=\min\ \{M,\nu(z)\}, \quad 
\nu^{[M]}_{\le k}(z):=\begin{cases}
	\nu^{[M]}(z)&\text{ if }\nu^{[M]}(z)\le k,\\
	0&\text{ if }\nu^{[M]}(z)> k.
\end{cases}
$$

Similarly, we define $\nu^{[M]}_{> k}$. We will write $N^{[M]}(r,\nu), \ N^{[M]}_{\le k} (r,\nu), \ N^{[M]}_{> k} (r,\nu)$ for $N(r, \nu^{[M]}),$ $N(r, \nu^{[M]}_{\le k}), \ N(r,\nu^{[M]}_{> k})$ as respectively. We will omit the character $^{[M]}$ if $M=\infty$.

\noindent
{\bf 2.3.}\ Let $f : \C^n \longrightarrow \P^N(\C)$ be a meromorphic mapping.
For arbitrarily fixed homogeneous coordinates
$(w_0 : \dots : w_N)$ on $\P^N(\C)$, we take a reduced representation
$f = (f_0 : \dots : f_N)$, which means that each $f_i$ is a  
holomorphic function on $\C^n$ and 
$f(z) = \big(f_0(z) : \dots : f_N(z)\big)$ outside the analytic set
$\{ f_0 = \dots = f_N= 0\}$ of codimension $\geq 2$. Let $a$ be a meromorphic mapping of $\C^n$ into $\P^N(\C)^*$ with reduced representation
$a = (a_0 : \dots : a_N)$. We denote by $T(r,f)$ the characteristic function of $f$ and by $m_{f,a}(r)$ the proximity function of $f$ with respect to $a$ (see \cite{Q}).

If  $(f,a)\not \equiv 0$, then the first main theorem for moving targets in value distribution theory states
$$T(r,f)+T(r,a)=m_{f,a}(r)+N_{(f,a)}(r).$$

\begin{theorem}[{The First Main Theorem for general position \cite[p. 326]{St1}}]
Let $f_j: \C^n\rightarrow \P^N(\C),\ 1\le j\le k$ be meromorphic mappings located in general position. Assume that $1\le k \le N.$ Then
$$N(r, \mu_{f_1\wedge\cdots\wedge f_{\lambda}})+m(r,f_1\wedge\cdots\wedge f_{\lambda})\le \sum_{1\le i\le\lambda}T(r,f_i)+O(1).$$
\end{theorem}

Here, by $\mu_{f_1\wedge\cdots\wedge f_{\lambda}}$ we denote the divisor associated with $f_1\wedge\cdots\wedge f_{\lambda}.$ We also denote by $N_{f_1\wedge\cdots\wedge f_{\lambda}}(r)$ the counting function associated with the divisor $\mu_{f_1\wedge\cdots\wedge f_{\lambda}}$.

Let $V$ be a complex vector space of dimension $N\ge 1.$ The vectors  $\{v_1,\cdots ,v_k\}$ are said to be in general position if for each selection of integers $1\le i_1<\cdots <i_p\le k$ with $p\le N,$ then  $v_{i_1}\wedge\cdots\wedge v_{i_p}\ne 0$. The vectors  $\{v_1,\cdots ,v_k\}$ are said to be in 
special position if they are not in general position. Take $1\le p \le k.$ Then $\{v_1,\cdots ,v_k\}$ are said to be in 
$p$-special position if for each selection of integers $1\le i_1<\cdots <i_p\le k,$ the vectors $v_{i_1},\cdots, v_{i_p}$ are in special position.

\begin{theorem}[{The Second Main Theorem for general position \cite[p. 320]{St1}}]
Let $M$ be a connected complex manifold of dimension $m.$ Let $A$ be a pure $(m-1)$-dimensional analytic subset of $M.$ Let $V$ be a complex vector space of dimension $n+1>1.$ Let $p$ and $k$ be integers with $1\le p\le k\le n+1.$ Let $f_j:M\rightarrow P(V),1\le j\le k,$ be meromorphic mappings. Assume that $f_1,...,f_k$ are in general position. Also assume that $f_1,...,f_k$ are in $p$-special position on $A.$ Then we have
$$\mu_{f_1\wedge\cdots\wedge f_k}\ge (k-p+1)\nu_{A}.$$
\end{theorem}
Here by $\nu_A$ we denote the reduced divisor whose support is the set A.

The following is a new second main theorem given by S. D. Quang \cite{Q}, which is an improvement the second main theorem of Thai-Quang in \cite{TQ}.
\begin{theorem}[{The Second Main Theorem for moving target \cite{Q}}]\label{2.3}
Let $f :\C^n \to \P^N(\C)$ be a meromorphic mapping. Let $\{a_i\}_{i=1}^q \ (q\ge 2N-m+2)$ be meromorphic mappings of $\C^n$ into $\P^N(\C)^*$ in general position such that $(f,a_i)\not\equiv 0\ (1\le i\le q),$ where $m+1=\mathrm{rank}_{\mathcal R\{a_i\}}(f)$. Then we have
$$|| \ \dfrac {q}{2N-m+2}T_f(r) \le \sum_{i=1}^q N_{(f,a_i)}^{[m]}(r) + o(T_f(r)) + O(\max_{1\le i \le q}T_{a_i}(r)).$$
\end{theorem}
 
As usual, by the notation ``$|| \ P$''  we mean the assertion $P$ holds for all $r \in [0,\infty)$ excluding a Borel subset $E$ of the interval $[0,\infty)$ with $\int_E dr<\infty$.

\vskip0.2cm 
\noindent
\textbf{Remark 2:} With the assumption of Theorem \ref{2.3}, we see that
\begin{align}\notag
||\ T_f(r) &\le \dfrac{m(2N-m+2)}{q}\sum_{i=1}^q N_{(f,a_i)}^{[1]}(r) + o(T_f(r)) + O(\max_{1\le i \le q}T_{a_i}(r))\\ 
\label{2.4}
& \le \dfrac{N(N+2)}{q}\sum_{i=1}^q N_{(f,a_i)}^{[1]}(r) + o(T_f(r)) + O(\max_{1\le i \le q}T_{a_i}(r)).
\end{align}

\section{Proofs of Main Theorems}

\noindent
{\bf 3.1. Proof of Theorem \ref{1.1}.}

It suffices to prove Theorem \ref{1.1} in the case of $\lambda \le N+1.$ 

Suppose that $f_1\wedge\cdots\wedge f_{\lambda}\not\equiv 0.$ We now prove the following.

\begin{claim}\label{Cl3.1}
For every $1\le i\le \lambda, \ 0\le j\le q-1$ and $1\le m\le N$,  we have
$$N^{[m]}_{\le k_j}(r,\nu^0_{(f_i,g_j)})\ge \dfrac{k_j+1}{k_j+1-m}N^{[m]}(r,\nu^0_{(f_i,g_j)})-\dfrac{m}{k_j+1-m}T(r,f_i).$$
\end{claim} 
Indeed, we have 
\begin{align*}
N^{[m]}_{\le k_j}(r,\nu_{(f_i,g_j)}^0)&=N^{[m]}(r,\nu_{(f_i,g_j)}^0)-N^{[m]}_{> k_j}(r,\nu_{(f_i,g_j)}^0)\\ 
&\ge N^{[m]}(r,\nu_{(f_i,g_j)}^0)-\dfrac{m}{k_j+1}N_{> k_j}(r,\nu_{(f_i,g_j)}^0)\\ 
&= N^{[m]}(r,\nu_{(f_i,g_j)}^0)-\dfrac{m}{k_j+1}N(r,\nu_{(f_i,g_j)}^0)+\dfrac{m}{k_j+1}N_{\le k_j}(r,\nu_{(f_i,g_j)}^0)\\
&\ge N^{[m]}(r,\nu_{(f_i,g_j)}^0)-\dfrac{m}{k_j+1}T(r,f_i)+\dfrac{m}{k_j+1}N^{[N]}_{\le k_j}(r,\nu_{(f_i,g_j)}^0).
\end{align*}
Thus
$$N^{[m]}_{\le k_j}(r,\nu_{(f_i,g_j)}^0)\ge \dfrac{k_j+1}{k_j+1-m}N^{[m]}(r,\nu^0_{(f_i,g_j)})-\dfrac{m}{k_j+1-m}T(r,f_i).$$
Then the claim is proved.

\begin{claim}\label{Cl3.2}
For every $1\le i\le \lambda,$  we have
$$\sum_{j=0}^{q-1}\min\{1,\nu^0_{(f_i,g_j),\le k_j}\}\le\dfrac{d}{\lambda -l+1}\mu_{f_1\wedge\cdots \wedge f_{\lambda}}(z)+q.\sum_{\beta}\mu_{g_{\beta (1)}\wedge\cdots\wedge g_{\beta (N+1)}}(z)$$
for each $z\not\in  A\bigcup_{i=1}^{\lambda}I(f_i),$ where the sum is over all injective maps $\beta: \{1,2,\cdots,N+1\}\to \{1,2,\cdots,q\}$
\end{claim}

Indeed, for each regular point $z_0\in\mathcal A\setminus (A\cup\bigcup_{i=1}^{\lambda}I(f_i)\cup\bigcup_{\beta\in T[N+1,q]}\{z|g_{\beta (1)}(z)\wedge\cdots\wedge g_{\beta (N+1)}(z)=0\})$ and for each increasing sequence  $1\le i_1<\cdots <i_l\le \lambda,$ we have 
$$f_{i_1}(z_0)\wedge\cdots\wedge f_{i_l}(z_0)=0.$$
By the Second Main Theorem for general position \cite[p. 320]{St1}, we have
$$\mu_{f_1\wedge\cdots\wedge f_{\lambda}}(z_0)\ge\lambda -(l-1).$$
Hence
\begin{align*}
\sum_{j=0}^{q-1}\min\{1,\nu^0_{(f_i,g_j),\le k_j}(z_0)\}\le \sum_{j=0}^{q-1}\min\{1,\nu^0_{(f_i,g_j)}(z_0)\}\le d \le\dfrac{d}{\lambda -l+1} \mu_{f_1\wedge\cdots\wedge f_{\lambda}}(z_0).
\end{align*}
If $z_0\in \bigcup_{\beta\in T[N+1,q]}\{z|g_{\beta (1)}(z)\wedge\cdots\wedge g_{\beta (N+1)}(z)=0\},$ then we have 
\begin{align*}
\sum_{j=0}^{q-1}\min\{1,\nu^0_{(f_i,g_j),\le k_j}(z_0)\}\le \sum_{j=0}^{q-1}\min\{1,\nu^0_{(f_i,g_j)}(z_0)\}\le q\sum_{\beta\in T[N+1,q]}\mu_{g_{\beta (1)}\wedge\cdots\wedge g_{\beta (N+1)}}(z_0).
\end{align*}
Thus, for each  $z\not\in A\cup\bigcup_{i=1}^{\lambda}I(f_i),$ we have
$$\sum_{j=0}^{q-1}\min\{1,\nu^0_{(f_i,g_j),\le k_j}(z)\}\le \dfrac{d}{\lambda -l+1} \mu_{f_1\wedge\cdots\wedge f_{\lambda}}(z)+q\sum_{\beta\in T[N+1,q]}\mu_{g_{\beta (1)}\wedge\cdots\wedge g_{\beta (N+1)}}(z).$$
The Claim \ref{Cl3.2} is proved.

The above Claim yields that
\begin{align*}
||\ \sum_{j=0}^{q-1}N^{[1]}_{\le k_j}(r,\nu^0_{(f_i,g_j)})&\le \dfrac{d}{\lambda -l+1} N_{f_1\wedge\cdots\wedge f_{\lambda}}(r)+q\sum_{\beta\in T[N+1,q]}N_{g_{\beta (1)}\wedge\cdots\wedge g_{\beta (N+1)}}(r)\\
&\le\dfrac{d}{\lambda -l+1}\sum_{i=1}^{\lambda}T(r,f_i)+q\sum_{\beta\in T[N+1,q]} \sum_{i=1}^{N+1}T(r,g_{\beta (i)})\\
&=\dfrac{d}{\lambda -l+1}T(r)+o(\max_{1\le i\le\lambda} T(r,f_i)).
\end{align*}
where $T(r)=\sum_{i=1}^\lambda T(r,f_i)$. Then, by Claim \ref{Cl3.1} and the inequality (\ref{2.4}) we have
\begin{align*}
||\ \dfrac{d\lambda}{\lambda -l+1}& T(r)\ge\sum_{i=1}^{\lambda}\sum_{j=0}^{q-1}N^{[1]}_{\le k_j}(r,\nu^0_{(f_i,g_j)})+o(T(r))\\
&\ge\sum_{i=1}^{\lambda}\sum_{j=0}^{q-1} \left(\dfrac{k_j+1}{k_j}N^{[1]}(r,\nu^0_{(f_i,g_j)})-\dfrac{1}{k_j}T(r,f_i)\right)+o(T(r))\\
&\ge\sum_{i=1}^{\lambda}\sum_{j=0}^{q-1}N^{[1]}(r,\nu^0_{(f_i,g_j)})-\sum_{j=0}^{q-1}\dfrac{1}{k_j}T(r)+o(T(r))\\
&\ge\sum_{i=1}^{\lambda}\dfrac{q}{N(N+2)}T(r,f_i)-\sum_{j=0}^{q-1}\dfrac{1}{k_j}T(r)+o(T(r))\\
&=\biggl (\dfrac{q}{N(N+2)}-\sum_{j=0}^{q-1}\dfrac{1}{k_j}\biggl)T(r)+o(T(r)).
\end{align*}
Letting $r\rightarrow +\infty,$ we get 
$$\sum\limits_{j=0}\limits^{q-1}\dfrac{1}{k_j} \ge \dfrac{q}{N(N+2)} - \dfrac{d\lambda}{\lambda-l+1}.$$
 This is a contradiction. Thus, $f_1\wedge\cdots\wedge f_{\lambda}=0.$ The theorem is proved.\hfill$\square$

\vskip0.3cm
\noindent
{\bf 3.2. Proof of Theorem \ref{1.3}.}

It suffices to prove Theorem \ref{1.3} in the case of $\lambda \le N+1.$  We set $m=\mathrm{rank}_{\mathcal R\{g_j\}}(f_i)-1$. Similar as \cite[Claim 3.1]{PP}, we have the following claim.
\begin{claim}\label{Cl3.3}
Let  $h_i: \C^n \to \P^N(\C)\ (1 \le i \le p \le N+1)$ be meromorphic mappings with reduced representations $h_i:=(h_{i0}:\cdots :h_{iN}).$ Let $a_i:\C^n \rightarrow \P^N(\C)^*\ (1\le i\le N+1)$ be moving hyperplanes with reduced representations $a_i:=(a_{i0}:\cdots :a_{iN}).$ 
Put $\tilde h_i:=((h_i,a_1):\cdots :(h_i,a_{N+1}))$
Assume that $a_1,\cdots, a_{N+1}$
are located in general position such that $(h_i,a_j)\not\equiv 0\ (1\le i \le p,\ 1\le j\le N+1).$ Let $S$ be a pure $(n-1)$-dimensional analytic subset of $\C^n$
such that $S\not\subset (a_1\wedge\cdots\wedge a_{N+1})^{-1}\{0\}.$ Then $h_1\wedge\cdots\wedge h_{p}=0$ on $S$ if and only if  
$\tilde h_1\wedge\cdots\wedge\tilde h_{p}=0$ on $S$.
\end{claim} 

We now prove Claim \ref{Cl3.3} .

\noindent
($\Rightarrow$)\ Suppose that $\{\tilde h_1\wedge\cdots\wedge\tilde h_{p}\}\not\equiv 0$ on $S.$ Then there exists $z_0\in S$ such that
$\tilde h_1(z_0)\wedge\cdots\wedge\tilde h _{p}(z_0)\ne 0 $. This means that the family $\{\tilde h_1(z_0),\cdots ,\tilde h_{p}(z_0)\}$ is linearly independent on $\C,$ i.e., the following matrix is of rank $p$
$$A=\left (\begin {array}{cccc}
(h_1,a_1)(z_0) &\cdots &(h_p ,a_{1})(z_0)\\
\vdots&\vdots&\vdots\\
(h_1,a_{N+1})(z_0) &\cdots &(h_p ,a_{N+1})(z_0)\\
\end {array}\right )$$
$$=\left (\begin {array}{cccc}
a_{10}(z_0) &\cdots &a_{1N}(z_0)\\
\vdots&\vdots&\vdots\\
a_{N+1 0}(z_0) &\cdots &a_{N+1 N}(z_0)\\
\end {array}\right )\cdot
\left (\begin {array}{cccc}
h_{10}(z_0) &\cdots &h_{p 0}(z_0)\\
\vdots&\vdots&\vdots\\
h_{1N}(z_0) &\cdots &h_{p N}(z_0)\\
\end {array}\right ).$$
Hence the matrix
$$\left (\begin {array}{cccc}
h_{10}(z_0) &\cdots &h_{p 0}(z_0)\\
\vdots&\vdots&\vdots\\
h_{1N}(z_0) &\cdots &h_{pN}(z_0)\\
\end {array}\right )$$
is of rank $p,$ i.e., $h_1(z_0)\wedge\cdots\wedge h_p(z_0)\ne 0.$ This yields that $h_1\wedge\cdots\wedge h_p\not\equiv 0$ on $S$. This is a contradiction.

\noindent
($\Leftarrow$) \ We see that the following matrix is of rank $\le p-1$ for each $z\in S$
$$A=\left (\begin {array}{cccc}
(h_1,a_1)(z) &\cdots &(h_p ,a_{1})(z)\\
\vdots&\vdots&\vdots\\
(h_1,a_{N+1})(z) &\cdots &(h_p ,a_{N+1})(z)\\
\end {array}\right ).$$
On the other hand, we have 
$$A=\left (\begin {array}{cccc}
a_{10}(z)&\cdots &a_{1N}(z)\\
\vdots&\vdots&\vdots\\
a_{N+1 0}(z) &\cdots &a_{N+1 N}(z)\\
\end {array}\right )\cdot
\left (\begin {array}{cccc}
h_{10}(z) &\cdots &h_{p 0}(z)\\
\vdots&\vdots&\vdots\\
h_{1N}(z) &\cdots &h_{pN}(z)\\
\end {array}\right ).$$
Since the family $\{a_i\}$ is located in general position and $S\not\subset (a_1\wedge\cdots\wedge a_{N+1})^{-1}\{0\},$ this implies that the matrix
$$\left (\begin {array}{cccc}
h_{10}(z) &\cdots &h_{p 0}(z)\\
\vdots&\vdots&\vdots\\
h_{1N}(z) &\cdots &h_{p N}(z)\\
\end {array}\right )$$
is of rank $\le p-1$ for each $z\in S, $ i.e.,  $h_1\wedge\cdots\wedge h_p\equiv 0$ on $S$.

The Claim \ref{Cl3.3} is proved.

\vskip0.2cm
We now continue to prove the theorem. Suppose that  $f_1\wedge\cdots\wedge f_{\lambda}\not\equiv 0$. 

For $\lambda$ indices $0= j_0<j_1<\cdots <j_{\lambda -1}\le N$ such that $\left (\begin {array}{cccc}
(f_1,g_{j_0}) &\cdots &(f_\lambda ,g_{j_0})\\
(f_1,g_{j_1}) &\cdots &(f_\lambda ,g_{j_1})\\
\vdots&\vdots&\vdots\\
(f_1,g_{j_{\lambda -1}}) &\cdots &(f_\lambda ,g_{j_{\lambda -1}})
\end {array}\right )$ is nondegenerate.

Put $J=\{j_0,\cdots,j_{\lambda-1}\} $, $J^c=\{0,\ldots ,q-1\}\setminus J$ and 
$B_J=\left (\begin {array}{cccc}
(f_1,g_{j_0}) &\cdots &(f_\lambda ,g_{j_0})\\
(f_1,g_{j_1}) &\cdots &(f_\lambda ,g_{j_1})\\
\vdots&\vdots&\vdots\\
(f_1,g_{j_{\lambda -1}}) &\cdots &(f_\lambda ,g_{j_{\lambda -1}})
\end {array}\right )$.

We now prove the following claim.

\begin{claim}\label{Cl3.4}
If $B_J$ is nondegenerate, i.e., $\det B_J\not\equiv 0$ then
\begin{align*}
\sum_{j\in J}(\min_{1\le i\le \lambda}\{\nu^0_{(f_i,g_j),\le k_j}\}&-\min\{1,\nu^0_{(f_1,g_j),\le k_j}\})\\
&+\sum_{j=0}^{q-1}(\lambda -l+1)\min\{1,\nu^0_{(f_1,g_j),\le k_j}\}\le \mu_{\tilde f_1\wedge\cdots\wedge \tilde f_{\lambda}}
\end{align*}
on the set $\C^n\setminus (A\cup\bigcup_{i=1}^{\lambda}I(f_i)\cup (g_{j_0}\wedge\cdots\wedge g_{j_{\lambda -1}})^{-1}(0)),$ where $\tilde f_i:=((f_i,g_{j_0}):\cdots :(f_i,g_{j_{\lambda-1}}))$ and $A=\bigcup_{0\le i<j\le q-1}Z_{(f_1,g_i)}\cap Z_{(f_1,g_j)}$.
\end{claim} 
Indeed,  put  $\mathcal A:=\bigcup_{j\in J}Z_{(f_1,g_j)}, \mathcal A^c:=\bigcup_{j\in J^c}Z_{(f_1,g_j)}.$ We consider the following two cases. 

{\bf Case 1.} \  Let $z_0\in \mathcal A\setminus (A\cup\bigcup_{i=1}^{\lambda}I(f_i)\cup(g_{j_0}\wedge\cdots\wedge g_{j_{\lambda-1}})^{-1}(0))$ be a regular point of $\mathcal A.$ Then $z_0$ is a zero of one of the meromorphic functions $\{(f_1,g_j)\}_{j\in J}.$ Without loss of generality we may assume that $z_0$ is a zero of $(f_1,g_{j_0}).$ Let $S$ be an irreducible component of $\mathcal A$ containing $z_0.$ Suppose that $U$ is an open neighborhood of $z_0$ in $\C^n$ such that $U\cap\{\mathcal A\setminus S\}=\emptyset$. Choose a holomorphic function $h$ on an open neighborhood $U'\subset U$ of $z_0$ such that  
$\nu_{h}(z)=\min_{1\le i\le \lambda}\{\nu^0_{(f_i,g_{j_0}),\le k_{j_0}}(z)\}$ if $z\in S$ and $\nu_{h}(z)=0$ if $z\not\in S.$ Then $(f_i,g_{j_0})=a_{i}h \ (1\le i\le\lambda),$ where  $a_{i}$ are holomorphic functions. Therefore, the matrix 
$\left (\begin {array}{cccc}
(f_1,g_{j_1}) &\cdots &(f_\lambda ,g_{j_1})\\
\vdots&\vdots&\vdots\\
(f_1,g_{j_{\lambda-1}}) &\cdots &(f_\lambda ,g_{j_{\lambda-1}})
\end {array}\right )$
is of rank $\le \lambda -1$. Hence, there exist $\lambda $ holomorphic functions $b_1,\cdots ,b_{\lambda}$, not all zeros, such that
$$\sum_{i=1}^{\lambda}b_i\cdot (f_i,g_{j_k})=0\ (1\le k\le\lambda-1).$$

Without loss of generality, we may assume that the set of common zeros of $\{b_i\}_{i=1}^{\lambda}$ is an analytic subset of codimension $\ge 2$. 
Then there exists an index $i_1,1\le i_1\le \lambda$, such that  $S\not\subset b_{i_1}^{-1}\{0\}$. We may assume that $i_1=\lambda.$ Then for each 
$z\in (U'\cup S)\setminus b_{\lambda}^{-1}\{0\},$ we have 
\begin{align*}
\tilde f_1(z)\wedge\cdots\wedge \tilde f_{\lambda}(z)&=\tilde f_1(z)\wedge\cdots\wedge \tilde f_{\lambda -1}(z)\wedge 
\biggl(\tilde f_{\lambda}(z)+\sum_{i=1}^{\lambda -1}\dfrac{b_i}{b_{\lambda}}\tilde f_i(z)\biggl)\\
&=\tilde f_1(z)\wedge\cdots\wedge \tilde f_{\lambda -1}(z)\wedge (V(z)h(z))\\
&=h(z)\cdot(\tilde f_1(z)\wedge\cdots\wedge \tilde f_{\lambda -1}(z)\wedge V(z)),
\end{align*}
where  $V(z):=(a_{\lambda}+\sum_{i=1}^{\lambda -1}\dfrac{b_i}{b_\lambda}a_i,0,\cdots ,0).$

By the assumption and by Claim \ref{Cl3.3}, for any increasing sequence $1\le i_1<\cdots <i_l\le \lambda -1,$ we have 
$\tilde f_{i_1}\wedge\cdots\wedge\tilde f_{i_l}=0$ on $S.$ This implies that the family $\{\tilde f_1,\cdots  \tilde f_{\lambda -1}, V\}$ is in $(l+1)$-special position on $S.$ By using The Second Main Theorem for general position \cite[p. 320]{St1}, we have
$$\mu_{\tilde f_1\wedge\cdots\wedge \tilde f_{\lambda -1}\wedge V}(z)\ge \lambda -l,\forall z\in  S.$$

Hence
\begin{equation*}
\mu_{\tilde f_1\wedge\cdots\wedge \tilde f_{\lambda}}(z)\ge\nu_h(z)+\lambda -l=\min_{1\le i\le \lambda}\{\nu^0_{(f_i,g_{j_0}),\le k_{j_0}}(z)\}+\lambda -l,
\end{equation*}
for all $z\in (U'\cup S)\setminus b_{i_1}^{-1}\{0\}$. This implies that
\begin{align*}
\sum_{j\in J}(\min_{1\le i\le \lambda}\{\nu^0_{(f_i,g_j),\le k_j}(z_0)\}&-\min\{1,\nu^0_{(f_1,g_j),\le k_j}(z_0)\})\\
&+\sum_{j=0}^{q-1}(\lambda -l+1)\min\{1,\nu^0_{(f_1,g_j),\le k_j}(z_0)\}\\
&=\min_{1\le i\le \lambda}\{\nu^0_{(f_i,g_{j_0}),\le k_{j_0}}(z_0)\}+\lambda -l\le \mu_{\tilde f_1\wedge\cdots\wedge \tilde f_{\lambda}}(z_0).
\end{align*}

{\bf Case 2.}\ Let $z_0\in \mathcal A^c\setminus (A\cup\bigcup_{i=1}^{\lambda}I(f_i)\cup(g_{j_0}\wedge\cdots\wedge g_{j_{\lambda-1}})^{-1}(0))$ be a regular point of $\mathcal A^c.$ Then  $z_0$ is a zero of $(f_1,g_j),j\in J^c.$ By the assumption and by Claim \ref{Cl3.3}, the family $\{\tilde f_1,\cdots ,\tilde f_{\lambda}\}$ is in $l$-special position on each irreducible component of $\mathcal A^c$ containing $z_0.$ By using The Second Main Theorem for general position \cite[p. 320]{St1}, we have 
$$\mu_{\tilde f_1\wedge\cdots\wedge \tilde f_{\lambda}}(z_0)\ge\lambda -l+1.$$
Hence
\begin{align*}
\sum_{j\in J}(\min_{1\le i\le \lambda}\{\nu^0_{(f_i,g_j),\le k_j}(z_0)\}&-\min\{1,\nu^0_{(f_1,g_j),\le k_j}(z_0)\})\\
&+\sum_{j=0}^{q-1}(\lambda -l+1)\min\{1,\nu^0_{(f_1,g_j),\le k_j}(z_0)\}\\
&=\lambda -l+1\le \mu_{\tilde f_1\wedge\cdots\wedge \tilde f_{\lambda}}(z_0).
\end{align*}

From the above two cases we get the desired inequality of the claim. 

\vskip0.2cm 
We now continue to prove the theorem. For each $j, 0\le j\le q-1,$ we set 
$$N_j(r)=\sum_{i=1}^{\lambda}N^{[m]}_{\le k_j}(r,\nu^0_{(f_i,g_j)})-((\lambda -1)m+1)N^{[1]}_{\le k_j}(r,\nu^0_{(f_1,g_j),}).$$
	
For each permutation $I=(j_0,\ldots ,j_{q-1})$ of $(0,\ldots ,q-1)$, we set 
$$T_I=\{r\in [1,+\infty ); N_{j_0}(r)\ge\cdots\ge N_{j_{q-1}}(r)\}.$$
It is clear that $\bigcup_{I}T_I=[1,+\infty)$. Therefore, there exists a permutation, for instance it is $I_0=(0,\ldots ,q-1)$, such that $\int_{T_{I_0}}dr =+\infty$. Then we have 
$$N_0(r)\ge N_1(r)\ge\cdots\ge N_{q-1}(r) \ \text{for all}\ r\in T_{I_0}.$$
	
By the assumption for $f_1\wedge \cdots\wedge f_\lambda\not\equiv 0$, there exist indices $J=\{j_0,\cdots, j_{\lambda-1}\}$ with $0=j_0< j_1< \cdots < j_{\lambda-1}\le N$ such that $\det B_J\not\equiv 0.$ We note that
 $$N_0(r)=N_{j_0}(r)\ge N_{j_1}(r)\ge\cdots\ge N_{j_{\lambda-1}}(r)\ge N_N(r), \ \text{for each} \ r\in T_{I_0}.$$

We see that $\min_{1\le i\le\lambda}a_i\ge\sum_{i=1}^{\lambda}\min\{m,a_i\}-(\lambda -1)m$
for every $\lambda$ non-negative integers $a_1,\ldots ,a_{\lambda}$. Then Claim \ref{Cl3.4} implies that
\begin{align*}
\sum_{j\in J}\Big(\sum_{i=1}^{\lambda}\min\{m,\nu^0_{(f_i,g_j),\le k_j}\}&-((\lambda -1)m+1)\min\{1,\nu^0_{(f_1,g_j),\le k_j}\}\Big)\\
&+\sum_{j=0}^{q-1}(\lambda -l+1)\min\{1,\nu^0_{(f_1,g_j),\le k_j}\}\le \mu_{\tilde f_1\wedge\cdots\wedge \tilde f_{\lambda}}.
\end{align*}
on the set $\C^n\setminus (A\cup\bigcup_{i=1}^{\lambda}I(f_i)\cup (g_{j_0}\wedge\cdots\wedge g_{j_{\lambda -1}})^{-1}(0)).$ 
Integrating both sides of this inequality, we have
\begin{align}\notag
\sum_{j\in J}\Big(\sum_{i=1}^{\lambda} N^{[m]}_{\le k_j}&(r,\nu^0_{(f_i,g_j)})-((\lambda -1)m+1)N^{[1]}_{\le k_j}(r,\nu^0_{(f_1,g_j)})\Big)\\
\label{3.1}
&+\sum_{j=0}^{q-1}(\lambda -l+1)N^{[1]}_{\le k_j}(r,\nu^0_{(f_1,g_j)})\le N_{\tilde f_1\wedge\cdots\wedge \tilde f_{\lambda}}(r) = N_{\det B_J}(r).
\end{align}
Also, by Jensen's formula, we have 
\begin{align}\label{3.2}
N_{\det B_J}(r)\le \int\limits_{S(r)} \log |\det B_J|\sigma_n + O(1)\le \sum_{i=1}^{\lambda}T(r,f_i)+o(\max_{1\le i\le \lambda}T(r,f_i)).
\end{align}
Set $T(r)=\sum_{i=1}^{\lambda}T(r,f_i)$. Combining (\ref{3.1}) and (\ref{3.2}), then for all $r\in I_0$, we have 
\begin{align*}
||\ T(r)&\ge\sum_{i=0}^{\lambda -1}N_{j_i}(r)+\sum_{j=0}^{q-1}(\lambda -l+1)N^{[1]}_{\le k_j}(r,\nu^0_{(f_1,g_j)})+o(T(r))\\
&\ge\dfrac{\lambda}{q}\sum_{j=0}^{q-1}N_j(r)+\sum_{j=0}^{q-1}(\lambda -l+1)N^{[1]}_{\le k_j}(r,\nu^0_{(f_1,g_j)})+o(T(r))\\
&=\sum_{j=0}^{q-1}\big(\lambda -l+1-\dfrac{\lambda ((\lambda -1)m+1)}{q}\big)N^{[1]}_{\le k_j}(r,\nu^0_{(f_1,g_j)})\\
&\ \ \ \ +\sum_{j=0}^{q-1}\dfrac{\lambda}{q}\sum_{i=1}^{\lambda}N^{[m]}_{\le k_j}(r,\nu^0_{(f_i,g_j)})+o(T(r))\\
&\ge\sum_{i=1}^{\lambda}\sum_{j=0}^{q-1}\big (\dfrac{\lambda}{q}+\dfrac{\lambda -l+1}{\lambda m}-\dfrac{(\lambda -1)m+1}{mq}\big )N^{[m]}_{\le k_j}(r,\nu^0_{(f_i,g_j)})+o(T(r))
\end{align*}

\begin{align*}
&\ge\sum_{i=1}^{\lambda}\sum_{j=0}^{q-1}\dfrac{q(\lambda -l+1)+\lambda (m-1)}{\lambda mq}\big(N^{[m]}(r,\nu^0_{(f_i,g_j)}) -\dfrac{m}{k_j+1-m}T(r,f_i)\big)+o(T(r))\\
&\ge\dfrac{q(\lambda -l+1)+\lambda (m-1)}{\lambda mq} \big (\sum_{i=1}^{\lambda}\sum_{j=0}^{q-1}N^{[m]}(r,\nu^0_{(f_i,g_j)})-\sum_{j=0}^{q-1}\dfrac{m}{k_j+1-m}T(r)\big )+o(T(r))\\
&\ge\dfrac{q(\lambda -l+1)+\lambda (m-1)}{\lambda mq}\big (\dfrac{q}{2N-m+2}-\sum_{j=0}^{q-1}\dfrac{m}{k_j+1-m}\big )T(r)+o(T(r)).
\end{align*}
Letting $r\longrightarrow +\infty,$ we have
$$ 1\ge \dfrac{q(\lambda -l+1)+\lambda (m-1)}{\lambda mq}\biggl (\dfrac{q}{2N-m+2}-\sum_{j=0}^{q-1}\dfrac{m}{k_j+1-m}\biggl ).$$
Thus
$$\sum_{j=0}^{q-1}\dfrac{1}{k_j+1-m}\ge \dfrac{q}{m(2N-m+2)}-\dfrac{\lambda q}{q(\lambda -l+1)+\lambda (m-1)}.$$
This is a contradiction. Thus, we have $f_1\wedge\cdots\wedge f_\lambda \equiv 0.$ \hfill$\square$

\end{document}